\documentclass[11pt]{amsart}
\RequirePackage{amsmath,amssymb,amsthm, amscd, comment,mathtools}

\usepackage{kantlipsum}
\usepackage[all]{xy}
\usepackage{graphicx}
\usepackage{enumitem}
\usepackage{times}
\usepackage{subfig}
\usepackage[pagebackref,breaklinks,colorlinks,linkcolor=red,anchorcolor=red,citecolor=blue]{hyperref}

\usepackage{enumitem}

\calclayout

\newtheorem{proposition}[subsubsection]{Proposition}
\newtheorem{corollary}[subsubsection]{Corollary}
\newtheorem{theorem}[subsubsection]{Theorem}
\newtheorem{lemma}[subsubsection]{Lemma}
\theoremstyle{definition}
\newtheorem{definition}[subsubsection]{Definition}
\theoremstyle{remark}
\newtheorem{recall}[subsubsection]{Recall}
\newtheorem{remark}[subsubsection]{Remark}

\numberwithin{equation}{subsubsection}

\DeclareMathAlphabet{\mathbbold}{U}{bbold}{m}{n}

\title{Algebraic $G$-theory in motivic homotopy categories}

\author{Fangzhou Jin}
\address{Fakult\"at f\"ur Mathematik\\
Universit\"at Duisburg-Essen\\
Thea-Leymann-Strasse 9\\
45127 Essen\\
Germany}
\email{\href{mailto:fangzhou.jin@uni-due.de}{fangzhou.jin@uni-due.de}}
\urladdr{\url{https://sites.google.com/site/fangzhoujin1/}}

\date{\number\day-\number\month-\number\year}


\begin{document}

\maketitle

\begin{abstract}

We prove that algebraic $G$-theory is representable in unstable and stable motivic homotopy categories; in the stable category we identify it with the Borel-Moore theory associated to algebraic $K$-theory, and show that such an identification is compatible with the functorialities defined by Quillen and Thomason.

\end{abstract}

\tableofcontents

\noindent

\section{Introduction}

\subsubsection{}
The study of Algebraic $G$-theory (or $K'$-theory) originated from Grothendieck's definition of the group $G_0$ in \cite[IV 2.2]{SGA6}, which associates to a noetherian scheme $X$ the Grothendieck group of the abelian category of coherent sheaves on $X$. It was Quillen who first defined higher $G$-groups for noetherian schemes in \cite[\S7]{Qui}, in the same way as he defined higher $K$-groups. He proved some properties of $G$-theory which fail to hold for $K$-theory in general, including homotopy invariance and localization property. This definition was generalized to all schemes by Thomason (\cite{TT}) by combining Waldhausen's construction (\cite{Wal}) and Grothendieck's work in \cite{SGA6}, which for noetherian schemes agrees with Quillen's definition. Thomason improved the functoriality of $G$-theory, and his definition is widely used together with his reformulation of algebraic $K$-theory.

\subsubsection{}
Ever since its birth, $G$-theory has been conceived as a cohomological theory that accompanies $K$-theory, and the two theories share many properties such as Brown-Gersten property and projective bundle formula. The definition of the two are quite similar as well: for example in Quillen's definition, $G$-theory is constructed in the same way as $K$-theory by replacing vector bundles by coherent sheaves; since any coherent sheaf over a regular noetherian scheme has a finite resolution by vector bundles (\cite[VIII 2.4]{SGA2}), it follows that Quillen $K$-theory and $G$-theory agree over regular schemes (which also turn out to be the same as Thomason $K$-theory).
We recall some important properties of (Thomason) $G$-theory, due to Quillen, Nisnevich and Thomason:
\begin{recall}[Properties of $G$-theory]
\label{G-property}

For any scheme $X$, denote by $G(X)$ Quillen-Thomason's $G$-theory spectrum (\cite[3.3]{TT}), whose homotopy groups $G_n(X)$ are $G$-theory groups. Then the following properties hold:
\begin{enumerate}

\item (contravariant functoriality) \label{G_pullback}  For any morphism of schemes $f:X\to X'$ of finite $\operatorname{Tor}$-dimension,
\footnote{The condition is to be understood in the sense of \cite[Section 7 2.5]{Qui}, where a morphism $f:X\to Y$ has finite $\operatorname{Tor}$-dimension if $\mathcal{O}_X$ has finite $\operatorname{Tor}$-dimension as a module over $f^{-1}(\mathcal{O}_Y)$, which is slightly stronger than the condition of being perfect, see \cite[III 4.1]{SGA6}.}
there is a map $f^*:G(X')\to G(X)$. (\cite[3.14.1]{TT})

\item (proper functoriality)
\label{G_proper_push}
  For any proper morphism between noetherian schemes $f:X\to X'$ there is a map $f_*:G(X)\to G(X')$. (\cite[3.16.1]{TT})

\item (homotopy invariance) 
\label{htp_G}
For any noetherian scheme $X$, the pull-back by the canonical projection induces a homotopy equivalence $G(X)\overset{p^*}{\to} G(X\times \mathbb{A}^1)$.  (\cite[\S7 4.1]{Qui})

\item (Brown-Gersten property)
\label{BG_G}
For any Cartesian diagram of noetherian schemes 
\begin{align}
\begin{gathered}
  \xymatrix@=10pt{
    V \ar[r]^-{} \ar[d] & Y \ar[d]^-{p}\\
    U \ar[r]^-{j} & X
  }
\end{gathered}
\end{align}
where $p$ is \'etale, $j$ an open immersion such that $p$ induces an isomorphism over $X-U$ (i.e. the square is a distinguished Nisnevich square), the square
\begin{align}
\begin{gathered}
  \xymatrix@=10pt{
    G(X) \ar[r]^-{j^*} \ar[d]_-{p^*} & G(U) \ar[d]^-{}\\
    G(Y) \ar[r]^-{} & G(V)
  }
\end{gathered}
\end{align}
is a homotopy pullback square. (\cite[4.4]{Nis})

\item (localization)
\label{localization_G}
Let $X$ be a noetherian scheme and $i:Z\to X$ a closed immersion with complementary open immersion $j:U\to X$. Then there is a homotopy fiber sequence
\begin{align}
G(Z)\overset{i_*}{\to}G(X)\overset{j^*}{\to}G(U).
\end{align}
In other words, there is a long exact sequence of $G$-theory groups
\begin{align}
\label{eq:G_loc}
\cdots\to G_{n+1}(U)\to G_{n}(Z)\overset{i_*}{\to}G_{n}(X)\overset{j^*}{\to}G_{n}(U)\cdots.
\end{align}
(\cite[\S7 3.2]{Qui})

\item (action by $K$-theory)
\label{action_K}
  For any noetherian scheme $X$ there is a map $K(X)\wedge G(X)\to G(X)$ which makes $G(X)$ a module over $K(X)$, which is compatible with the contravariant functoriality. (\cite[3.15.3]{TT})

\item (coincidence with $K$-theory for regular schemes)
\label{agree_regular}
  For any regular noetherian scheme $X$, there is a canonical homotopy equivalence $K(X)\simeq G(X)$. (\cite[3.21]{TT})

\item (projective bundle formula) 
\label{proj_bundle_G}
Let $X$ be a noetherian scheme and $\mathcal{E}$ be a vector bundle of rank $r$ over $X$. Denote by $p:\mathbb{P}(\mathcal{E})\to X$ the projection of the projective bundle associated to $E$ and $\alpha\in K_0(\mathbb{P}((\mathcal{E}))$ the class of the line bundle $\mathcal{O}(1)$. 
\footnote{Throughout the article, we use the convention $\mathbb{P}(\mathcal{E})=\operatorname{Proj}(\operatorname{Sym}\mathcal{(\mathcal{E}}^\vee))$, the projectivization of the symmetric algebra of the dual of its sheaf of sections.}
Then the map
\begin{align}
\begin{split}
        \bigoplus_{i=0}^{r-1} G_n(X)&\to G_n(\mathbb{P}(E))\\
        (x_0,\cdots,x_{r-1})&\mapsto \sum_{i=0}^{r-1}\alpha^i\cdot p^*(x_i)
\end{split}
\end{align}
is an isomorphism. (\cite[\S7 4.3]{Qui})
\end{enumerate}
\end{recall}

\subsubsection{}
On the other hand, in the framework of motivic homotopy theory developed by Morel and Voevodsky (\cite{MV}), it is known that \emph{Borel-Moore theories} have similar functorialities: they are covariant with respect to proper morphisms and contravariant with respect to local complete intersection morphisms, and possess a localization sequence in the form of~\eqref{eq:G_loc} 
 (see \cite{Deg2}, \cite{Jin}). Since $G$-theory is a module over $K$-theory ring spectrum and agrees with the latter for regular schemes, the properties enlisted in Recall~\ref{G-property} suggest that $G$-theory should be considered as a Borel-Moore theory associated to $K$-theory seen as a cohomology theory. The main purpose of this paper is to prove such a result, 
 using the six functors formalism established in \cite{Ayo} and \cite{CD}.

\subsubsection{}
The first problem concerning $G$-theory in motivic homotopy categories is representability problem, which is already well-studied for Thomason's $K$-theory: Morel and Voevodsky showed in \cite[Theorem 4.3.13]{MV} that over a regular noetherian base scheme $S$, $K$-theory is represented in the pointed $\mathbb{A}^1$-homotopy category $\mathbf{H}_\bullet(S)$ by the group completion of the infinite Grassmannian space $\mathbb{Z}\times Gr$; the representability can be promoted to the stable homotopy category $\mathbf{SH}(S)$, where $K$-theory is represented by a $\mathbb{P}^1$-spectrum $\mathbf{KGL}_S$ constructed from $\mathbb{Z}\times Gr$ in \cite{Voe3} (see \cite{Rio}, \cite{PPR}).

\subsubsection{}
In Section~\ref{chap_rep_G} we follow the same steps to deal with $G$-theory over an arbitrary noetherian base scheme $S$: we first show the representability of $G$-theory in $\mathbf{H}_\bullet(S)$ (Corollary~\ref{rep_H}). The model we construct uses axioms of a module over the $K$-theory spectrum of free vector bundles, which is used to construct a strict model for the proper functoriality. We then imitate Voevodsky's construction to obtain an object $\mathbf{GGL}_S$ in $\mathbf{SH}(S)$ 
(Definition~\ref{bott_class}), show that $\mathbf{GGL}_S$ is indeed an $\Omega$-spectrum (Corollary~\ref{stabilization_iso}), and therefore deduce the representability of $G$-theory in $\mathbf{SH}(S)$ by $\mathbf{GGL}_S$:
\begin{theorem}[see Corollary~\ref{rep_SH}]
\label{intro_rep_SH}
The spectrum $\mathbf{GGL}_S$ represents the $G$ theory in the stable homotopy category $\mathbf{SH}(S)$. In other words, for any smooth $S$-scheme $X$, there is a natural isomorphism
\begin{align}
Hom_{\mathbf{SH}(S)}(\Sigma^\infty X_+[n], \mathbf{GGL}_S)\simeq G_n(X).
\end{align}
\end{theorem}

\subsubsection{}
Then it comes to study the functorialities of $G$-theory as a Borel-Moore theory. Apart from the proper functoriality, there is a refined Gysin functoriality for abstract Borel-Moore theories, which has been constructed in \cite{Deg2} in for oriented theories (see Recall~\ref{recall_BM}) and recently generalized in \cite{DJK} to a more general setting. We would like to identify these functorialities with the ones of $G$-theory in Recall~\ref{G-property}.

\subsubsection{}
In Section~\ref{chap_func}, we study functorial behaviors of the spectrum $\mathbf{GGL}_S$. In Section~\ref{contra_func} we establish its contravariance with respect to the exceptional inverse image functor $f^!$. There are two simple cases: on the one hand, it follows from the proper functoriality of $G$-theory and our definition that for any proper morphism $p:W\to X$, there is a canonical map $\phi_p:p_*\mathbf{GGL}_W\to \mathbf{GGL}_X$; on the other hand, for any open immersion $j:U\to X$ there is a canonical identification $j^*\mathbf{GGL}_X\simeq \mathbf{GGL}_U$ by restriction to an open subset. We show that these two cases can be glued, and with some extra work studying the case of projective morphisms, we obtain:
\begin{theorem}[see Proposition~\ref{exceptional_G}]
\label{intro_functorality_G}
For any separated morphism of finite type $f:Y\to X$, there is a canonical isomorphism $\mathbf{GGL}_Y\simeq f^!\mathbf{GGL}_X$.

\end{theorem}

\subsubsection{}
In Section~\ref{section:sm_func}, we show that for smooth morphisms, the map obtained above agrees with another isomorphism obtained from Bott periodicity and purity isomorphisms (see Proposition~\ref{sm_embed}). The proof uses compatibility between Gysin morphisms to show the quasi-projective case, and deduce the general case by noetherian induction. In Section~\ref{appl} we apply the results above to show the compatibility between functorialities: over a regular base scheme, we identify $G$-theory with the Borel-Moore theory associated to $K$-theory, and deduce the following result:
\begin{theorem}[see Corollary~\ref{G=BM_func}]
\label{BM_G}
For $S$ a regular scheme, $f:X\to S$ a separated morphism of finite type and integers $m,n$, we have an isomorphism $\mathbf{KGL}^{BM}_{n,m}(X/S)\simeq G_{n-2m}(X)$, compatible with proper covariance and lci 
\footnote{We say that a morphism of schemes $f:X\to Y$ is \emph{local complete intersection} (abbreviated as ``lci'') if it factors as the composition of a regular closed immersion followed by a smooth morphism.}
 contravariance on both sides.

\end{theorem}

\subsubsection{}
Note that the result above has been used in \cite{Deg2} to establish a Riemann-Roch theorem for singular schemes, which generalizes \cite[Theorem 18.3]{Ful}, see Remark~\ref{remark_RR} below.

\subsubsection*{\bf Acknowledgments}
The author would like to thank Fr\'ed\'eric D\'eglise suggesting this problem during his PhD thesis, and Denis-Charles Cisinski for many discussions and ideas. He would like to thank Alberto Navarro for a remark on orientations, and Adeel Khan for pointing out some errors in a preliminary version.

\subsubsection*{\bf Notations and Conventions}
\begin{enumerate}
\item Throughout the article we assume that all schemes are noetherian. 
\item A \emph{smooth} morphism stands for a separated smooth morphism of finite type. We denote by $Sch$ the category of schemes, and for any scheme $X$, we denote by $Sch_X$ (resp. $Sm_X$) the category of separated $X$-schemes of finite type (resp. the category of smooth $X$-schemes). 
\item For any pair of adjoint functors $(F,G)$ between two categories, we denote by $ad_{(F,G)}:1\to GF$ and $ad'_{(F,G)}:GF\to 1$ the unit and conuit maps of the adjunction.
\end{enumerate}

\section{Representability of algebraic $G$-theory}
\label{chap_rep_G}
In this section we prove the representability of algebraic $G$-theory in motivic homotopy categories. We construct a special model which will be used later to study the functorial properties.

\subsection{Construction of the model}

\subsubsection{}
In what follows we deal with the action by algebraic $K$-theory. Following \cite{TT}, the action of the $K$-theory spectrum over $G$-theory stems from the tensor product of a cohomologically bounded pseudo-coherent complex by a perfect complex. We would like to rectify this action into a morphism of strict presheaves. To do this, we start with the rectification of the tensor product by a free $\mathcal{O}_X$-module. This is done by strictifying the direct sum operation. 

\subsubsection{}
For any scheme $X$, we denote by $Coh^{fl}(X)$ the category of cohomologically bounded, pseudo-coherent, bounded below complexes of flasque quasi-coherent $\mathcal{O}_X$-modules (\cite[3.11.5]{TT}). If $S$ is a scheme, any $S$-morphism between smooth $S$-schemes is lci (\cite[B.7.6]{Ful}), therefore of finite Tor-dimension (\cite[VIII 1.7]{SGA6}).

\begin{definition}
\label{def_C}
Let $S$ be a scheme. For any scheme $X\in Sm_S$, we denote by $Sm_S/X$ the full subcategory of $Sm_S$ whose objects are also $X$-schemes, i.e. the overcategory of $X$-objects in $Sm_S$. We define $C_S(X)$ as the category where
\begin{itemize}
\item Objects are data of
\begin{itemize}
\item for every $Y\in Sm_S/X$, a complex $E_Y$ in $Coh^{fl}(Y)$,
\item for every morphism $p:Y\to Z$ in $Sm_S/X$, a map of complexes $\alpha_{p}:E_Z\to p_{*}E_Y$, 
\end{itemize}
such that
\begin{enumerate}
\item for every morphism $p_{}:Y\to Z$ in $Sm_S/X$, the map $\tilde{\alpha}_{p}:Lp^*_{}E_Z\to E_Y$ obtained from $\alpha_{p}$ by adjunction in the derived category is a quasi-isomorphism,
\item for every $Y\in Sm_S/X$, the map $\alpha_{Id_{Y,Y}}:E_Y\to E_Y$ is the identity map,
\item for composable morphisms 
\begin{align}
Y
\xrightarrow{p}
Z
\xrightarrow{q}
W
\end{align}
in $Sm_S/X$, the composition
\begin{align}
E_W
\xrightarrow{\alpha_q}
q_*E_Z
\xrightarrow{q_*\alpha_p}
q_*p_*E_Y
\simeq
(q\circ p)_*E_Y
\end{align}
agrees with the map $\alpha_{q\circ p}$.
\end{enumerate}
\item A morphism $\eta$ from $((E_Y)_{Y\in Sm_S/X},(\alpha_p)_{p:Y\to Z\in Sm_S/X})$ to $((F_Y)_{Y\in Sm_S/X},(\beta_p)_{p:Y\to Z\in Sm_S/X})$ is the data of a map of complexes $\eta_Y:E_Y\to F_Y$ for every $Y\in Sm_S/X$, such that for every morphism $p:Y\to Z$ in $Sm_S/X$, the following diagram commutes:
\begin{align}
\begin{gathered}
  \xymatrix@=14pt{
 E_Z \ar[r]^-{\eta_Z} \ar[d]_-{\alpha_p} & F_Z \ar[d]^-{\beta_p} \\
 p_{*}E_Y \ar[r]^-{p_{*}\eta_Y} & p_{*}F_Y.
  }
\end{gathered}
\end{align}
\end{itemize}

\end{definition}

\subsubsection{}
According to the definition in \cite[1.2.11]{TT}, the category $C_S(X)$ is a complicial biWaldhausen category with respect to the abelian category of all diagrams of the form $A_Z\to p_{*}A_Y$ where $p:Y\to Z$ is a morphism in $Sm_S/X$ and $A_Y$ (resp. $A_Z$) is a $\mathcal{O}_Y$-module (resp. $\mathcal{O}_Z$-module). 
Recall the following result in \cite[Th\'eor\`eme 2.15]{Cis}:
\begin{proposition}
\label{prop:cis2.14}
Let $F:\mathcal{C}\to\mathcal{D}$ be a right exact functor between Waldhausen categories (\cite[1.2.5]{TT}) such that
\begin{enumerate}
\item Every morphism $u$ in $\mathcal{C}$ and $\mathcal{D}$ factors as $u=pi$ where $i$ is a cofibration and $p$ is a weak equivalence.
\item Every morphism in $\mathcal{C}$ whose image in $Ho(\mathcal{C})$ is an isomorphism is indeed a weak equivalence, and the same property holds for $\mathcal{D}$.
\item $F$ induces an equivalence between the homotopy categories $Ho(\mathcal{C})\simeq Ho(\mathcal{D})$.
\end{enumerate}
Then the induced map on $K$-theory spaces $K(F):K(\mathcal{C})\to K(\mathcal{D})$ is a homotopy equivalence.
\end{proposition}

\begin{lemma}
\label{C_coh}
For any scheme $X$, the canonical map between Waldhausen categories
\begin{align}
\begin{split}
        a_S(X):C_S(X)&\to Coh^{fl}(X)\\
        ((E_Y)_{Y\in Sm_S/X},(\alpha_p)_{p:Y\to Z\in Sm_S/X})&\mapsto E_X
\end{split}
\end{align}
induces a homotopy equivalence between the associated Waldhausen $K$-theory spaces (\cite[1.5.2]{TT}). Consequently, by \cite[3.11]{TT}, for any scheme $X$ the $K$-theory space associated to $C(X)$ is homotopy equivalent to the algebraic $G$-theory space (\cite[3.3]{TT}).

\end{lemma}

\proof

It suffices to check that the conditions in Proposition~\ref{prop:cis2.14} are satisfied. The first two axioms are standard (see \cite[Exemple 2.5]{Cis}). For the third axiom, we construct a functor $b_S(X):Coh^{fl}(X)\to C_S(X)$ which is an inverse of $a_S(X)$ up to quasi-isomorphism: for $E\in Coh^{fl}(X)$ let $P(E)$ be a functorial resolution of $E$ by a cohomologically bounded, pseudo-coherent, bounded above complex of flat quasi-coherent $\mathcal{O}_X$-modules, and let $b_S(X)$ send $E$ to the following data:
\begin{itemize}
\item for every $Y\in Sm_S/X$, the complex $Go(g^*P(E))$ in $Coh^{fl}(Y)$, where $Go$ is the Godement resolution functor, 
\item for every morphism $p:Y\to Z$ in $Sm_S/X$, the canonical map 
\begin{align}
Go(g^*P(E))
\to
p_*Go(p^*g^*P(E))
\simeq
p_*Go((g\circ p)^*P(E)).
\end{align}
\end{itemize}
The functor $b_S(X)$ is well-defined since any morphism $f$ in $Sm_S$ has finite $\operatorname{Tor}$-dimension, and therefore $f^*$ preserves cohomological boundedness and pseudo-coherence. It is easy to see that $b_S(X)$ is an inverse to $a_S(X)$ up to quasi-isomorphism, and the result follows.
\endproof

\begin{lemma}
\label{func_CX}
\begin{enumerate}
\item The map $X\mapsto C_S(X)$ defines a (strict) presheaf of sets over the category $Sm_S$.
\item Let $f:T\to S$ be a proper morphism. For any $S$-scheme $X$, denote by $f_X:X_T=X\times_ST\to X$ the base change of $f$. Then there is a covariant map $f_{X*}:C_T(X_T)\to C_S(X)$, such that for any morphism $p:Y\to X$ in $Sm_S$, the following diagram commutes:
\begin{align}
\label{diag:C_pf_pb}
\begin{gathered}
  \xymatrix@=10pt{
 C_T(Y_T)  \ar[d]_-{f_{Y*}} & C_T(X_T) \ar[l]_-{p_T^*} \ar[d]^-{f_{X*}} \\
 C_S(Y)  & C_S(X) \ar[l]_-{p^*}
  }
\end{gathered}
\end{align}
where $p_T:X_T\to Y_T$ is the base change of $p$.
\end{enumerate}
\end{lemma}
\proof
\begin{enumerate}
\item Let $q:W\to X$ be a morphism in $Sm_S$. 
The functor $Lq^*$ preserves cohomologically bounded pseudo-coherent complexes, and the composition with $q$ identifies $Sm_S/W$ as a full subcategory of $Sm_S/X$. Then there is a well-defined map $q^*:C_S(X)\to C_S(W)$ which sends the data $((E_Y)_{Y\in Sm_S/X},(\alpha_p)_{p:Y\to Z\in Sm_S/X})$ to the subdata $((E_Y)_{Y\in Sm_S/W},(\alpha_p)_{p:Y\to Z\in Sm_S/W})$. 
Such a map defines a strict contravariant functoriality of the map $X\mapsto C_S(X)$ on $Sm_S$.

\item 
We define the map $f_{X*}:C_T(X_T)\to C_S(X)$ by sending the data $((E_Y)_{Y\in Sm_T/{X_T}},(\alpha_p)_{p:Y\to Z\in Sm_T/{X_T}})$ to the following data:
\begin{itemize}
\item for every $Y\in Sm_S/X$, the complex of $\mathcal{O}_Y$-modules $f_{Y*}E_{Y_T}$,
\item for every morphism $p:Y\to Z$ in $Sm_S/X$, the map 
\begin{align}
f_*\alpha_{p_T}:
f_{Z*}E_{Z_T}
\xrightarrow{f_{Z*}\alpha_{p_T}}
f_{Z*}p_{T*}E_{Y_T}
\simeq
p_{*}f_{Y*}E_{Y_T}.
\end{align}
\end{itemize}
Since $f_Y:Y_T\to Y$ is a proper morphism between noetherian schemes, the functor $Rf_{Y*}$ on the derived category of $\mathcal{O}_{Y_T}$-modules preserves flasqueness, pseudo-coherence and cohomological boundedness (\cite[3.16]{TT}). 
The map $Lp^*f_{Z*}E_{Z_T}\to f_{Y*}E_{Y_T}$ agrees with the composition
\begin{align}
\label{eq:pf_bc}
Lp^*f_{Z*}E_{Z_T}
=
Lp^*Rf_{Z*}E_{Z_T}
\to
Rf_{Y*}Lp_T^*E_{Z_T}
\xrightarrow{Rf_{Y*}\tilde{\alpha}_{p_T}}
Rf_{Y*}E_{Y_T}
=
f_{Y*}E_{Y_T}.
\end{align}
The map $Lp^*Rf_{Z*}E_{Z_T}\to Rf_{Y*}Lp_T^*E_{Z_T}$ is a quasi-isomorphism by the $\operatorname{Tor}$-independent base change theorem (\cite[2.5.6]{TT}), and by assumptions the map~\eqref{eq:pf_bc} is a quasi-isomorphism. Therefore the map $f_{X*}:C_T(X_T)\to C_S(X)$ is well-defined. The commutativity of the diagram~\eqref{diag:C_pf_pb} follows directly from the construction.
\end{enumerate}
\endproof

\subsubsection{}
\label{section:Cplus}
By \cite[XI.3. Theorem 1]{ML}, every monoidal category $C$ is equivalent to a strict monoidal category $C^\oplus$ in a strongly monoidal way. For any scheme $S$ and $X\in Sm_S$ we denote by $C_S(X)^\oplus$ the strict monoidal category corresponding to the monoidal category $(C_S(X),\oplus)$. In other words, in the category $C_S(X)^\oplus$ direct sums are strictly associative, i.e. for any triple $(A,B,C)$ of objects we have an identification $(A\oplus B)\oplus C=A\oplus(B\oplus C)$. This also makes $C_S(X)^\oplus$ a Waldhausen category in a canonical way. By Lemma~\ref{C_coh}, we have the following:
\begin{lemma}
\label{lem:C=G}
For any scheme $X$, there is a canonical homotopy equivalence between the Waldhausen $K$-theory space of $C_S(X)^\oplus$ and Thomason's $G$-theory space $G(X)$.

\end{lemma}

\subsubsection{}
For any scheme $X$, denote by $Vect^{f}_X$ the category of free vector bundles (of finite rank) over $X$ as a subcategory of $Vect_X$. The category $Vect^{f}_X$ is a Waldhausen category by taking isomorphisms as weak equivalences and monomorphisms as cofibrations, and has a monoidal structure by tensor products. We define an action of $Vect^{f}_X$ on the category $C_S(X)^\oplus$ by setting
\begin{align}
\label{eq:action_free+}
\begin{split}
        Vect^{f}_X\times C_S(X)^\oplus&\to C_S(X)^\oplus\\
        (V, E)&\mapsto \underbrace{E\oplus E\cdots E\oplus E}_\text{$rank(V)$ times}
\end{split}
\end{align}
The map~\eqref{eq:action_free+} induces a map between $K$-theory spectra
\begin{equation}
\label{action_K_free}
K(Vect^{f}_X)\wedge K(C_S(X)^\oplus)\to K(C_S(X)^\oplus).
\end{equation}

\subsection{Representability in motivic homotopy categories}
\subsubsection{}
In what follows, we fix a noetherian scheme $S_0$ and work over schemes in $Sch^{}_{S_0}$. 
Denote by $F$ the functor in \cite[Lemma 2.2]{RSO} over the category $Sch^{}_{S_0}$, which is a lax symmetric monoidal fibrant replacement functor on the category of pointed motivic spaces over $S$ for any $S\in Sch^{}_{S_0}$.
\begin{definition}
Let $S$ be a scheme of finite type over $S_0$. Denote by $K^0_S$ (respectively $K^{f}_S$) the pointed presheaf $X\mapsto K(Vect_X)$ (respectively $X\mapsto K(Vect^{f}_X)$) over the category $Sm_S$, which has a monoidal structure induced by tensor product in $Vect_X$ (respectively $Vect^{f}_X$). Denote by $K_S=F(K^0_S)$, which has the structure of a monoid. 
\end{definition}

\subsubsection{}
The canonical map $K^{f}_S\to K^0_S$ is a Nisnevich local weak equivalence in the category of pointed motivic spaces over $Sm/S$, since every vector bundle is Zariski locally free. It follows that the canonical map $K^f_S\xrightarrow{\sim} K_S$ is also a weak equivalence. Since this map is compatible with the structure of monoids of both objects, by \cite[Theorem 4.3]{SS} we have the following:
\begin{lemma}
\label{equiv_K_free}
The derived functors of restriction and base change induce equivalences between homotopy categories
\begin{equation}
Ho(K^{f}_S-Mod)\simeq Ho(K_S-Mod).
\end{equation}
\end{lemma}

\subsubsection{}
\label{G_unstable}
Denote by $G^\oplus_S$ the presheaf $X\mapsto K(C_S(X)^\oplus)$ over the category $Sm_S$. The map~\eqref{action_K_free} endows $G^\oplus_S$ with the structure of a module over $K^{f}_S$. By Lemma~\ref{equiv_K_free}, there is a cofibrant replacement $\bar{G}^\oplus_S$ of $G^\oplus_S$ in $K^{f}_S-Mod$ which is weakly equivalent to the base change $G^\oplus_S\otimes_{K^{f}_S}K_S$. 
The object $G^\oplus_S\otimes_{K^{f}_S}K_S$ is a module over $K_S$, and we choose $G_S$ as a fibrant replacement of the object $\bar{G}^\oplus_S$ in the category $K_S-Mod$.

\subsubsection{}
The object $G_S$ in the unstable motivic homotopy category $\mathbf{H}(S)$ represents Thomason-Trobaugh $G$-theory:
\begin{lemma}
\label{rep_H}
For any $X\in Sm_S$ there is a canonical isomorphism:
\begin{equation}
\label{rep_G_H}
[S^n\wedge X_+, G_S]_{\mathbf{H}(S)}\simeq G_n(X).
\end{equation}
\end{lemma}

\proof

Since $G$-theory satisfies Nisnevich descent and homotopy invariance, by Lemma~\ref{lem:C=G} the object $G^\oplus_S$ represents $G$-theory in $\mathbf{H}(S)$, that is, for any $X\in Sm/S$ there is a functorial isomorphism:
\begin{equation}
[S^n\wedge X_+, G^\oplus_S]_{\mathbf{H}(S)}\simeq G_n(X).
\end{equation}
Then the lemma follows from the fact that the two objects $G^\oplus_\cdot$ and $G_\cdot$ are isomorphic in the homotopy category $\mathbf{H}_\bullet(S)$.
\endproof

\begin{lemma}
\label{pf_G}
For any proper morphism $f:W\to X$ in $Sch^{}_{S_0}$, there is a canonical map $\phi_f:f_*G_W\to G_X$ in the unstable motivic homotopy category $\mathbf{H}(S)$. If $g:V\to W$ is another proper morphism, then the composition
\begin{equation}
(f\circ g)_*G_V\simeq f_*g_*G_V\xrightarrow{f_*\phi_g}f_*G_W\xrightarrow{\phi_f} G_X
\end{equation}
agrees with $\phi_{(f\circ g)}$.
\end{lemma}
\proof

Since the construction in~\ref{section:Cplus} is functorial, for any proper morphism $f:W\to X$ and $Y\in Sm_X$, the map $f_*:C_W(Y_W)\to C_X(Y)$ in Lemma~\ref{func_CX} induces a map $f_*:C_W(Y_W)^\oplus\to C_X(Y)^\oplus$ which is an exact functor between Waldhausen categories. For any morphism $p:Y\to Z$ in $Sm_X$ with a Cartesian diagram
\begin{align}
\begin{gathered}
  \xymatrix@=10pt{
 Y_W \ar[r]^-{p_W} \ar[d]_-{f_Y} & Z_W \ar[r]^-{} \ar[d]^-{f_Z} & W \ar[d]^-{f} \\
 Y \ar[r]^-{p} & Z \ar[r]^-{} & X
  }
\end{gathered}
\end{align} 
the following diagram commutes by Lemma~\ref{func_CX}:
\begin{align}
\begin{gathered}
  \xymatrix@=10pt{
 C_W(Y_W)^\oplus  \ar[d]_-{f_{Y*}} & C_W(Z_W)^\oplus \ar[l]_-{p_W^*} \ar[d]^-{f_{Z*}} \\
 C_X(Y)^\oplus  & C_X(Z)^\oplus. \ar[l]_-{p^*}
  }
\end{gathered}
\end{align} 
Applying the $K$-theory space functor, we get a functorial map of presheaves of $K$-theory spaces $f_*G^\oplus_W\to G^\oplus_X$. By construction the two objects $G^\oplus_\cdot$ and $G_\cdot$ are isomorphic in the homotopy category $\mathbf{H}(S)$, and the result follows.
\endproof

\subsubsection{}
By definition, $K^0_S$ is the presheaf of Quillen $K$-theory spectra, and since Quillen $K$-theory agrees with Thomason $K$-theory for affine schemes (\cite[3.9]{TT}), there is Nisnevich local weak equivalence between $K^0_S$ and the presheaf of Thomason $K$-theory. 
Therefore for any $X\in Sm_S$, there is a canonical map
\begin{equation}
\label{K_to_KS}
K_n(X)\to [S^n\wedge X_+,K_S]_{\mathbf{H}_\bullet(S)}
\end{equation}
which is an isomorphism if $S$ is regular (see also \cite[Theorem 6.5]{Voe3}). Now for any vector bundle $\mathcal{E}$ of rank $r$ over $S$, let $x$ be the class $[\mathcal{O}(1)]$ in $K_0(\mathbb{P}(\mathcal{E}\oplus\mathcal{O}_X))$, and denote by $\nu(\mathcal{E})$ the class
\begin{equation}
\label{gen_bott_class}
\nu(\mathcal{E})=x^r-[\wedge^1\mathcal{E}]x^{r-1}+\cdots+(-1)^r[\wedge^r\mathcal{E}]\in K_0(\mathbb{P}(\mathcal{E}\oplus\mathcal{O}_X)).
\end{equation}
By abuse of notation, we still denote by $\nu(\mathcal{E})$ its image via the map~\eqref{K_to_KS}. By \cite[Proposition 3.2.17]{MV} we have an isomorphism
\begin{equation}
\label{Thom_proj_bundle}
Th(\mathcal{E})
\simeq
\mathbb{P}(\mathcal{E}\oplus\mathcal{O}_X)/\mathbb{P}(\mathcal{E}).
\end{equation}
Since the restriction of $\nu(\mathcal{E})$ to $\mathbb{P}(\mathcal{E})$ is zero,  $\nu(\mathcal{E})$ induces a class in $[Th(\mathcal{E}),K_S]_{\mathbf{H}_\bullet(S)}$. Via the canonical maps $K_S\wedge K_S\to K_S$ and $K_S\wedge G_S\to G_S$, the class of $\nu(\mathcal{E})$ induces maps in $\mathbf{H}_\bullet(S)$:
\begin{equation}
\label{bott_K}
K_S\to\underline{Hom}(Th(\mathcal{E}),K_S),
\end{equation}
\begin{equation}
\label{bott_G}
G_S\to\underline{Hom}(Th(\mathcal{E}),G_S).
\end{equation}

\begin{lemma}
\label{omega}
The maps~\eqref{bott_K} and~\eqref{bott_G} are isomorphisms.

\end{lemma}
\proof

The case of $K_S$ is known by \cite[Lemma 6.1.3.3]{Rio2}. The case of $G_S$ is similar: for any $X\in Sm/S$, the map~\eqref{bott_G} induces a map
\begin{equation}
\label{omega_spec_map}
G_n(X)
\overset{\eqref{rep_G_H}}{\simeq}
[S^n\wedge X_+,G_S]_{\mathbf{H}_\bullet(S)}
\to
[Th(\mathcal{E})\wedge S^n\wedge X_+,G_S]_{\mathbf{H}_\bullet(S)},
\end{equation}
where by the isomorphism~\eqref{Thom_proj_bundle}, the group $[Th(\mathcal{E})\wedge S^n\wedge X_+,G_S]_{\mathbf{H}_\bullet(S)}$ is identified with the $n$-th homotopy group of the homotopy fiber of the canonical map
\begin{equation}
G(\mathbb{P}(\mathcal{E}\oplus\mathcal{O}_S)\times_SX)\to G(\mathbb{P}(\mathcal{E})\times_SX).
\end{equation}
By definition, the composition
\begin{equation}
G_n(X)
\xrightarrow{\eqref{omega_spec_map}}
[Th(\mathcal{E})\wedge S^n\wedge X_+,G_S]_{\mathbf{H}_\bullet(S)}
\to
G_n(\mathbb{P}(\mathcal{E}\oplus\mathcal{O}_S)\times_SX)
\end{equation}
is given by $a\mapsto \nu\cdot p^*a$. Therefore the map~\eqref{omega_spec_map} is an isomorphism by the projective bundle formula for $G$-theory (Recall~\ref{G-property} \eqref{proj_bundle_G}) and the result follows.
\endproof

\subsubsection{}
In particular, if $\mathcal{E}$ is a trivial bundle of rank $1$, we have an isomorphism $Th(\mathbb{A}^1)\simeq\mathbb{P}^1$ as pointed motivic spaces, and $\nu(\mathbb{A}^1_S):\mathbb{P}^1\to K_S$ is the image of the Bott class $[\mathcal{O}(1)]-1\in K_0(\mathbb{P}^1_S)$. Following \cite[D\'efinition IV.1]{Rio} we stabilize the $G$-theory sheaf as follows:
\begin{definition}
\label{bott_class}
For any $S\in Sch^{}_{S_0}$, we define $\mathbf{KGL}_S$ (respectively $\mathbf{GGL}_S$) as the $\mathbb{P}^1$-spectrum with $\mathbf{KGL}_{S,n}=K_S$ (respectively $\mathbf{GGL}_{S,n}=G_S$) for all $n$ and suspension maps
\begin{equation}
\mathbb{P}^1\wedge K_S\xrightarrow{\nu(\mathbb{A}^1_S)\wedge1}K_S\wedge K_S\xrightarrow{} K_S
\end{equation}
\begin{equation}
\mathbb{P}^1\wedge G_S\xrightarrow{\nu(\mathbb{A}^1_S)\wedge1}K_S\wedge G_S\xrightarrow{} G_S.
\end{equation}
\end{definition}

\subsubsection{}
Both spectra $\mathbf{KGL}_S$ and $\mathbf{GGL}_S$ are fibrant for the projective model structure since they are degreewise fibrant. We now determine the cohomology theories they represent in the stable homotopy category. The case of $\mathbf{KGL}_S$ has already been well studied in the literature: by \cite[Theorem 3.6]{Voe3}, the spectrum $\mathbf{KGL}_S$ is isomorphic to the Voevodsky spectrum (\cite[6.2]{Voe3}) in $\mathbf{SH}(S)$, as explained in \cite[Remark 1.2.2]{PPR}. By \cite[Theorem 6.9]{Voe3} we have the following:
\begin{corollary}
\label{rep_KGL}
For any smooth $X\in Sm/S$, there is a functorial isomorphism
\begin{equation}
[\Sigma^\infty X_+[n], \mathbf{KGL}_S]_{\mathbf{SH}(S)}\simeq KH_n(X)
\end{equation}
where the right hand side is the homotopy $K$-theory group (\cite[Definition IV.12.7]{Wei})
\end{corollary}

We now go back to the spectrum $\mathbf{GGL}_S$. As a consequence of Lemma~\ref{omega} we have the following:
\begin{corollary}
\label{stabilization_iso}
The spectrum $\mathbf{GGL}_S$ is an $\Omega$-spectrum in the stable homotopy category $\mathbf{SH}(S)$. Consequently, for any pointed motivic space $A$ over $S$ the canonical map
\begin{equation}
[A,G_S]_{\mathbf{H}_\bullet(S)}\to [\Sigma^{\infty}A,\mathbf{GGL}_S]_{\mathbf{SH}(S)}
\end{equation}
is an isomorphism.
\end{corollary}
From Lemma~\ref{rep_H} and Corollary~\ref{stabilization_iso} we deduce that $\mathbf{GGL}$ represents algebraic $G$-theory in $\mathbf{SH}$:
\begin{corollary}
\label{rep_SH}
The spectrum $\mathbf{GGL}_S$ represents the $G$-theory in the stable homotopy category $\mathbf{SH}(S)$: for any $X\in Sm/S$, there is a canonical isomorphism
\begin{equation}
\label{GGL_G}
[\Sigma^\infty X_+[n], \mathbf{GGL}_S]_{\mathbf{SH}(S)}\simeq G_n(X).
\end{equation}

\end{corollary}

Since Thom spaces are invertible for the $\wedge$-product in $\mathbf{SH}$, by Lemma~\ref{omega} we have the following:
\begin{corollary}[Bott periodicity]
\label{Bott_per}
For any vector bundle $\mathcal{E}$ over $S$, there is a canonical isomorphism in $\mathbf{SH}(S)$:
\begin{equation}
\label{Bott_per_map}
\mathbf{GGL}_S\wedge Th(\mathcal{E})\simeq\mathbf{GGL}_S.
\end{equation}
\end{corollary}

\section{Functorial properties of the $G$-theory spectrum}
\label{chap_func}
In this section we study functorial properties of the $G$-theory spectrum, especially its contravariance with respect to the exceptional inverse image functor $f^!$. Using such a result we identify algebraic $G$-theory with the Borel-Moore theory associated to algebraic $K$-theory.

\subsection{Contravariant functoriality}
\label{contra_func}
\begin{lemma}
\label{GGL_open}
For any smooth morphism $f:Y\to X$, there is a canonical isomorphism 
\begin{equation}
\label{sm_func}
\chi_f:\mathbf{GGL}_Y\simeq f^*\mathbf{GGL}_X.
\end{equation}
\end{lemma}
\proof
For any $W\in Sm/Y$, by Corollary~\ref{rep_SH} we have a canonical isomorphism
\begin{equation}
[\Sigma^\infty_YW_+[n], \mathbf{GGL}_Y]_{\mathbf{SH}(Y)}\simeq G_n(W).
\end{equation}
Since $f$ is smooth, the functor $f^*$ has a left adjoint $f_\#$ such that $f_\#\Sigma^\infty_YW_+=\Sigma^\infty_XW_+$ (\cite[Proposition 3.2.9]{MV}), and therefore we have
\begin{align}
\begin{split}
[\Sigma^\infty_YW_+[n], f^*\mathbf{GGL}_X]_{\mathbf{SH}(Y)}
&=
[f_\#\Sigma^\infty_YW_+[n], \mathbf{GGL}_X]_{\mathbf{SH}(X)}\\
&=
[f_\#\Sigma^\infty_YW_+[n], \mathbf{GGL}_X]_{\mathbf{SH}(X)}\\
&\simeq
[\Sigma^\infty_XW_+[n], \mathbf{GGL}_X]_{\mathbf{SH}(X)}
=
G_n(W),
\end{split}
\end{align}
and the result follows.
\endproof

\subsubsection{}
The following lemma is straightforward from Lemma~\ref{pf_G} by stabilization:
\begin{lemma}
\label{pf_GGL}
For any proper morphism $f:W\to X$, there is a map $\phi_f:f_*\mathbf{GGL}_W\to \mathbf{GGL}_X$. If $g:V\to W$ is another proper morphism, then the composition
\begin{align}
(f\circ g)_*\mathbf{GGL}_V
\simeq
f_*g_*\mathbf{GGL}_V
\xrightarrow{f_*\phi_g}
f_*\mathbf{GGL}_W
\xrightarrow{\phi_f}
\mathbf{GGL}_X
\end{align}
agrees with $\phi_{(f\circ g)}$.
\end{lemma}

\subsubsection{}
The map $\phi_f$ can be understood as follows: for any $V\in Sm/X$, the map on $G$-groups
\begin{equation}
\begin{split}
G_n(W\times_XV)
&\simeq
[\Sigma^\infty_W(W\times_XV)_+[n], \mathbf{GGL}_W]_{\mathbf{SH}(W)}
=
[f^*\Sigma^\infty_XV_+[n], \mathbf{GGL}_W]_{\mathbf{SH}(W)}\\
&=
[\Sigma^\infty_XV_+[n], f_*\mathbf{GGL}_W]_{\mathbf{SH}(X)}
\xrightarrow{\phi_f}
[\Sigma^\infty_XV_+[n], \mathbf{GGL}_X]_{\mathbf{SH}(X)}
\simeq
G_n(V)
\end{split}
\end{equation}
induced by $\phi_f$ agrees with the proper functoriality of $G$-theory (Recall~\ref{G-property} \eqref{G_proper_push}). The following lemma follows from the construction:
\begin{lemma}
\label{comp_prop_sm}
For any Cartesian square of schemes
\begin{align}
\begin{split}
  \xymatrix@=10pt{
    V \ar[r]^-{q} \ar[d]_-{g} & Y \ar[d]^-{f}\\
    W \ar[r]^-{p} & X
  }
\end{split}
\end{align}
where $f$ and $g$ are smooth and $p$ and $q$ are proper, the following diagram commutes:
\begin{align}
\begin{split}
  \xymatrix@=10pt{
 q_*g^*\mathbf{GGL}_W  & q_*\mathbf{GGL}_V \ar[r]^-{\phi_q} \ar[l]_-{\chi_g}^-{\sim} & \mathbf{GGL}_Y \ar[d]^-{\chi_f}_-{\wr}\\
 f^*p_*\mathbf{GGL}_W \ar[rr]^-{\phi_p} \ar[u]^-{\wr} &  & f^*\mathbf{GGL}_X.
  }
\end{split}
\end{align}
\end{lemma}

\subsubsection{}
\label{not_com_GGL}
We use the following notation: given two composable morphisms $Z\xrightarrow{g}Y\xrightarrow{f}X$ and two maps $a:\mathbf{GGL}_Y\to f^!\mathbf{GGL}_X$, $b:\mathbf{GGL}_Z\to g^!\mathbf{GGL}_Y$, we denote by $a\cdot b$ the composition
\begin{align}
\mathbf{GGL}_Z\xrightarrow{b} g^!\mathbf{GGL}_Y\xrightarrow{g^!a}g^!f^!\mathbf{GGL}_X\simeq (f\circ g)^!\mathbf{GGL}_X.
\end{align}

\subsubsection{}
If $j:U\to X$ is an open immersion, by Lemma~\ref{GGL_open} there is a canonical isomorphism
\begin{equation}
\label{oi_func}
\chi_j:\mathbf{GGL}_U\simeq j^*\mathbf{GGL}_X=j^!\mathbf{GGL}_X.
\end{equation}
On the other hand, by Lemma~\ref{pf_GGL} and adjunction, for any proper morphism $f:W\to X$, we have a map
\begin{equation}
\label{prop_func}
\psi_f:\mathbf{GGL}_W\to f^!\mathbf{GGL}_X.
\end{equation}
It is then clear that the formation of the maps~\eqref{oi_func} and~\eqref{prop_func} are compatible with compositions of open immersions and proper morphisms respectively. These two types of maps are compatible by the following lemma:
\begin{lemma}
\label{comp_oi_prop}
For any commutative square of schemes
\begin{align}
\begin{split}
  \xymatrix@=10pt{
    V \ar[r]^-{q} \ar[d]_-{k} & U \ar[d]^-{j}\\
    Y \ar[r]^-{p} & X
  }
\end{split}
\end{align}
where $j$ and $k$ are open immersions and $p$ and $q$ are proper, we have $\chi_j\cdot\psi_q=\psi_p\cdot\chi_k$.
\end{lemma}

\proof

If the diagram is Cartesian, the result follows from Lemma~\ref{comp_prop_sm}. In the general case the canonical morphism $V\to Y\times_XU$ is an open and closed immersion, and the result follows from the Cartesian case.
\endproof

\subsubsection{}
Since every separated morphism of finite type $f$ has a compactification by our assumptions (\cite{Con}), namely a factorization $f=p\circ j$ where $p$ is proper and $j$ is an open immersion, using the technique in \cite[XVII]{SGA4}, we have the following functoriality by gluing maps~\eqref{oi_func} and~\eqref{prop_func}:
\begin{proposition}
\label{psi_comp}
There is a unique family of maps $\psi_f:\mathbf{GGL}_Y\to f^!\mathbf{GGL}_X$ associated to all separated morphisms of finite type $f:Y\to X$ such that
\begin{enumerate}
\item \label{oi_agree}
If $f$ is an open immersion, $\psi_f=\chi_f$ is the map \eqref{oi_func};
\item \label{prop_agree}
If $f$ is proper, $\psi_f$ is the map \eqref{prop_func};
\item \label{comp_stable}
For any two composable morphisms $f\circ g=h$, we have $\psi_f\cdot\psi_g=\psi_h$.
\end{enumerate}
\end{proposition}

\proof

We know that for every morphism the category of its compactifications is non empty and left filtering. For a morphism $f$ with $f=p\circ j$ a compactification, we set $\psi_f=\psi_p\cdot\xi_j$. By Lemma~\ref{comp_oi_prop}, the map $\psi_f$ is well-defined and independent on the choice of compactification. The properties~\eqref{oi_agree} and~\eqref{prop_agree} and the uniqueness are satisfied by definition. The property~\eqref{comp_stable} follows by applying Lemma~\ref{comp_oi_prop} again.
\endproof

\begin{proposition}
\label{exceptional_G}
For any separated morphism of finite type $f$, the map $\psi_f$ is an isomorphism.

\end{proposition}

The statement is local, and by localizing $f$, we only need to show the case where $f$ is quasi-projective. Therefore by Proposition~\ref{psi_comp} we only need to deal with three cases: open immersions, closed immersions, and the projection of a projective space. The open immersion case is Lemma~\ref{GGL_open}, and the two remaining cases will follow from Proposition~\ref{ci_psi} and Proposition~\ref{psi_proj} below.

\begin{proposition}
\label{ci_psi}
Let $X$ be a scheme and $i:Z\to X$ a closed immersion with complementary open immersion $j:U\to X$. Then 
\begin{enumerate}
\item \label{localization_htpfiber}
The map $\phi_i:i_*\mathbf{GGL}_Z\to\mathbf{GGL}_X$ identifies $i_*\mathbf{GGL}_Z$ canonically with homotopy fiber of the canonical map $\mathbf{GGL}_X\to j_*\mathbf{GGL}_U$;
\item The map $\psi_i:\mathbf{GGL}_Z\to i^!\mathbf{GGL}_X$ is an isomorphism.
\end{enumerate}
\end{proposition}

\proof
\begin{enumerate}
\item For any smooth $X$-scheme $X'$, denote $Z'=Z\times_XX'$ and $U'=U\times_XX'$. By Corollary~\ref{rep_SH} we have
\begin{align}
[\Sigma^\infty X'_+[n], i_*\mathbf{GGL}_Z]_{\mathbf{SH}(X)}=[i^*\Sigma^\infty X'_+[n], \mathbf{GGL}_Z]_{\mathbf{SH}(Z)}\simeq G_n(Z'),
\end{align}
\begin{align}
[\Sigma^\infty X'_+[n], j_*\mathbf{GGL}_U]_{\mathbf{SH}(X)}=[j^*\Sigma^\infty X'_+[n], \mathbf{GGL}_U]_{\mathbf{SH}(U)}\simeq G_n(U').
\end{align}
Then the result follows from the localization sequence (Recall~\ref{G-property} \eqref{localization_G}).
\item Using the localization distinguished triangle
\begin{equation}
\label{loc_triangle}
i_!i^!\to1\to j_*j^*\xrightarrow{}i_!i^![1]
\end{equation}
(\cite[Proposition 1.4.9]{Ayo}), by Lemma~\ref{GGL_open} we know that the map $i_*\mathbf{GGL}_Z\to i_!i^!\mathbf{GGL}_X$ induced by $\psi_i$ is an isomorphism. The result then follows from the fact that the functor $i_*=i_!$ is conservative (\cite[1.4.1]{Ayo}).
\end{enumerate}
\endproof

\begin{definition}
\label{def_xi}
For any smooth morphism $f:Y\to X$, by Lemma~\ref{GGL_open}, Corollary~\ref{Bott_per} and purity for $f$, we denote by $\xi_f:\mathbf{GGL}_Y\simeq f^!\mathbf{GGL}_X$ the following composition of isomorphisms
\begin{equation}
\label{eq:xi_sm}
\xi_f:
\mathbf{GGL}_Y
\overset{\eqref{Bott_per_map}}{\simeq}
\mathbf{GGL}_Y\wedge Th(T_f)
\overset{\chi_f}{\simeq}
f^*\mathbf{GGL}_X\wedge Th(T_f)
\overset{}{\simeq}
f^!\mathbf{GGL}_X.
\end{equation}
\end{definition}

\begin{proposition}
\label{psi_proj}
Let $p:\mathbb{P}^r_X\to X$ be the projection of a projective space. Then $\psi_p:\mathbf{GGL}_{\mathbb{P}^r_X}\to p^!\mathbf{GGL}_X$ is an isomorphism.

\end{proposition}

\proof
It suffices to show that $\xi_p=\psi_p$ where $\xi_p:\mathbf{GGL}_{\mathbb{P}^r_X}\simeq p^!\mathbf{GGL}_X$ is as defined in~\eqref{eq:xi_sm}. By considering the diagram
$$
\resizebox{\textwidth}{!}{%
  \xymatrix{
    \mathbf{GGL}_{\mathbb{P}^r_X} \ar @{-}[r]^-{\eqref{Bott_per_map}}_-{\sim} \ar[dd]_-{ad_{(p_!,p^!)}} & \mathbf{GGL}_{\mathbb{P}^r_X}\wedge Th(T_p) \ar[d]_-{ad_{(p_!,p^!)}} \ar[r]^-{\chi_p}_-{\sim} & p^*\mathbf{GGL}_X\wedge Th(T_p) \ar[r]^-{ad_{(p_!,p^!)}} \ar@/^.5cm/[rdd]^-{\sim} & p^!p_*(p^*\mathbf{GGL}_X\wedge Th(T_p)) \ar @{-}[d]^-{\wr}  \\
      & p^!p_*(\mathbf{GGL}_{\mathbb{P}^r_X}\wedge Th(T_p)) \ar @{-}[r]^-{\sim} \ar[rru]^-{\chi_p}_-{\sim} & p^!p_\#\mathbf{GGL}_{\mathbb{P}^r_X} \ar[r]^-{\chi_p}_-{\sim} & p^!p_\#p^*\mathbf{GGL}_X \ar[d]^-{ad'_{(p_\#,p^*)}}\\
 p^!p_*\mathbf{GGL}_{\mathbb{P}^r_X} \ar @{-}[ru]^-{\eqref{Bott_per_map}}_-{\sim} \ar[rrr]^-{\phi_p} &  &  & p^!\mathbf{GGL}_X
  }%
}
$$
we see that it suffices to show that the composition
\begin{align}
\begin{split}
p_*\mathbf{GGL}_{\mathbb{P}^r_X}
&\overset{\eqref{Bott_per_map}}{\simeq}
p_*(\mathbf{GGL}_{\mathbb{P}^r_X}\wedge Th(T_p))
\overset{}{\simeq}
p_\#\mathbf{GGL}_{\mathbb{P}^r_X}\\
&\overset{\chi_p}{\simeq}
p_\#p^*\mathbf{GGL}_X
\xrightarrow{ad'_{(p_\#,p^*)}}
\mathbf{GGL}_X
\end{split}
\end{align}
agrees with $\phi_p$. Again by considering the diagram
$$
  \xymatrix{
    p_\#\mathbf{GGL}_{\mathbb{P}^r_X} \ar@/^-1.7cm/[dd]_-{\wr} \ar[rr]^-{\chi_p}_-{\sim} \ar[d]^-{ad_{(p_\#,p^*)}} &  & p_\#p^*\mathbf{GGL}_X \ar@/^1.5cm/[dd]^-{ad'_{(p_\#,p^*)}}\\
    p_\#p^*p_\#\mathbf{GGL}_{\mathbb{P}^r_X} \ar @{-}[r]^-{\sim} & p_\#p^*(p_*\mathbf{GGL}_{\mathbb{P}^r_X}\wedge Th(T_p)) \ar[ld]_-{ad'_{(p_\#,p^*)}} \ar @{-}[r]^-{\eqref{Bott_per_map}}_-{\sim} & p_\#p^*p_*\mathbf{GGL}_{\mathbb{P}^r_X} \ar[ld]_-{ad'_{(p_\#,p^*)}} \ar[u]^-{\phi_p}\\
    p_*\mathbf{GGL}_{\mathbb{P}^r_X}\wedge Th(T_p) \ar @{-}[r]^-{\eqref{Bott_per_map}}_-{\sim} & p_*\mathbf{GGL}_{\mathbb{P}^r_X} \ar[r]^-{\phi_p} & \mathbf{GGL}_X
  }
$$
we see that it suffices to show that the composition
\begin{equation}
\begin{split}
\label{GGL_Pr}
\mathbf{GGL}_{\mathbb{P}^r_X}
&\xrightarrow{ad'_{(p_\#,p^*)}}
p^*p_\#\mathbf{GGL}_{\mathbb{P}^r_X}
\overset{}{\simeq}
p^*p_*(\mathbf{GGL}_{\mathbb{P}^r_X}\wedge Th(T_p))\\
&\overset{\eqref{Bott_per_map}}{\simeq}
p^*p_*\mathbf{GGL}_{\mathbb{P}^r_X}
\overset{\chi_p}{\simeq}
p^*p_*p^*\mathbf{GGL}_X
\xrightarrow{ad'_{(p^*,p_*)}}
p^*\mathbf{GGL}_X
\end{split}
\end{equation}
agrees with the map $\chi_p$. For any scheme $Y$ smooth over $\mathbb{P}^r_X$, we have canonical isomorphisms
\begin{align}
[\Sigma^\infty_{\mathbb{P}^r_X}Y_+[n],\mathbf{GGL}_{\mathbb{P}^r_X}]_{\mathbf{SH}(\mathbb{P}^r_X)}
\simeq
G_n(Y)
\simeq
[\Sigma^\infty_{\mathbb{P}^r_X}Y_+[n],p^*\mathbf{GGL}_X]_{\mathbf{SH}(\mathbb{P}^r_X)},
\end{align}
\begin{align}
\begin{split}
&[\Sigma^\infty_{\mathbb{P}^r_X}Y_+[n],p^*p_*\mathbf{GGL}_{\mathbb{P}^r_X}]_{\mathbf{SH}(\mathbb{P}^r_X)}
=
[p_\#\Sigma^\infty_{\mathbb{P}^r_X}Y_+[n],p_*\mathbf{GGL}_{\mathbb{P}^r_X}]_{\mathbf{SH}(\mathbb{P}^r_X)}\\
=
&[\Sigma^\infty_XY_+[n],p_*\mathbf{GGL}_{\mathbb{P}^r_X}]_{\mathbf{SH}(X)}
=
[p^*\Sigma^\infty_XY_+[n],\mathbf{GGL}_{\mathbb{P}^r_X}]_{\mathbf{SH}(\mathbb{P}^r_X)}\\
=
&[\Sigma^\infty_{\mathbb{P}^r_X}\mathbb{P}^r_{Y+}[n],\mathbf{GGL}_{\mathbb{P}^r_X}]_{\mathbf{SH}(\mathbb{P}^r_X)}
\simeq
G_n(\mathbb{P}^r_Y).
\end{split}
\end{align}
Therefore applying the functor $[\Sigma^\infty_{\mathbb{P}^r_X}Y_+[n],\cdot]_{\mathbf{SH}(\mathbb{P}^r_X)}$, the map~\eqref{GGL_Pr} becomes
\begin{equation}
\label{G_pbf}
G_n(Y)\xrightarrow{p^*_Y}G_n(\mathbb{P}^r_Y)\xrightarrow{p_{Y*}}G_n(Y)
\end{equation}
where $p_Y:\mathbb{P}^r_Y\to Y$ is the canonical projection. By projection formula for $G$-theory (\cite[\S7 2.10]{Qui}) and the description of the proper push-forward of a projective space on $K_0$ (\cite[VI 5.2]{SGA6}), for any $a\in G_n(Y)$ we have
\begin{align}
p_{Y*}p^*_Ya=p_{Y*}1\cdot a=a
\end{align}
where $1\in K_0(\mathbb{P}^r_Y)$ is the unit element. Therefore the map~\eqref{G_pbf} is indeed the identity map, and by Lemma~\ref{GGL_open}, the map~\ref{GGL_Pr} agrees with the map $\chi_p$, which finishes the proof.
\endproof
This finishes the proof of Proposition~\ref{exceptional_G}, which can also be stated in the form of Theorem~\ref{intro_functorality_G}.

\subsection{On the smooth functoriality}
\label{section:sm_func}
In this section we show that the smooth functoriality in Definition~\ref{def_xi} and the one in Proposition~\ref{psi_comp} are compatible.
\subsubsection{}
Recall that by Definition~\ref{def_xi}, for any smooth morphism $f:X\to S$, there is a canonical isomorphism $\xi_f:\mathbf{GGL}_Y\simeq f^!\mathbf{GGL}_X$. In this section we check its compatibility with the map $\psi_f$. We use the notations in~\ref{not_com_GGL} for compositions of maps between $\mathbf{GGL}$.
\begin{lemma}
\label{xi_comp}
The map $p\mapsto\xi_p$ is compatible with composition of smooth morphisms.
\end{lemma}
\proof
The statement follows from the fact that purity isomorphisms, Bott periodicity maps and the map $\chi_p$ are compatible with compositions.
\endproof

\begin{lemma}
\label{comp_ci_sm}
For any Cartesian square
\begin{align}
\begin{split}
  \xymatrix@=10pt{
    Y \ar[r]^-{k} \ar[d]_-{g} \ar@{}[rd]|{\Delta} & X \ar[d]^-{f}\\
    T \ar[r]_-{i} & S
  }
\end{split}
\end{align}
where $i$ and $k$ are closed immersions and $f$ and $g$ are smooth, then $\xi_f\cdot\psi_k=\psi_i\cdot\xi_g$. 
\end{lemma}

\proof
The result follows from Lemma~\ref{comp_oi_prop} and the fact that purity isomorphisms and Bott periodicity isomorphisms are compatible with base change.
\endproof

\begin{lemma}
\label{xi_section}
If $f:X\to S$ be a smooth morphism with $i:S\to X$ a section, then we have $\xi_f\cdot\psi_i=1$. In other words the composition
\begin{align}
\mathbf{GGL}_S
\overset{\psi_i}{\simeq}
i^!\mathbf{GGL}_X
\overset{\xi_f}{\simeq}
i^!f^!\mathbf{GGL}_S
\simeq
\mathbf{GGL}_S
\end{align}
is the identity map.
\end{lemma}
\proof

By considering the diagram
\begin{equation}
\begin{split}
  \xymatrix@=13pt{
     & i^!i_*\mathbf{GGL}_S \ar[r]^-{\phi_i} \ar[d]^{\xi_f}_-{\wr} & i^!\mathbf{GGL}_X \ar[d]^{\xi_f}_-{\wr}\\
    \mathbf{GGL}_S \ar[ru]^-{ad_{(i_*,i^!)}} \ar@/^-.5cm/[rr]^-{\sim} & i^!i_*i^!f^!\mathbf{GGL}_S \ar[r]^-{ad'_{(i_*,i^!)}} & i^!f^!\mathbf{GGL}_S
  }
\end{split}
\end{equation}
we are reduced to show that the commutativity of the following diagram:
\begin{equation}
\begin{split}
\label{xi_section_custom}
  \xymatrix@=13pt{
     i_*i^!f^!\mathbf{GGL}_S \ar[r]^-{ad'_{(i_*,i^!)}} \ar[d]^{}_-{\wr} & f^!\mathbf{GGL}_S \\
    i_*\mathbf{GGL}_S \ar[r]^-{\phi_i} & \mathbf{GGL}_X. \ar[u]^{\xi_f}_-{\wr}
  }
\end{split}
\end{equation}
Denote by $j:U\to X$ the open complement to $i:S\to X$. By Proposition~\ref{ci_psi}, the map $\phi_i$ identifies $i_*\mathbf{GGL}_S$ with the homotopy fiber of the canonical map $\mathbf{GGL}_X\to j_*\mathbf{GGL}_U$. By localization triangle~\eqref{loc_triangle}, the upper horizontal map in the diagram~\eqref{xi_section_custom} identifies $i_*i^!f^!\mathbf{GGL}_S$ with the homotopy fiber of the canonical map $f^!\mathbf{GGL}_S\xrightarrow{ad_{(j^!,j_*)}}j_*j^!f^!\mathbf{GGL}_S\overset{(\xi_{f\circ j})^{-1}}{\simeq}j_*\mathbf{GGL}_U$. Therefore the left vertical map in the diagram~\eqref{xi_section_custom} is a canonical isomorphism between the two homotopy fibers $i_*i^!f^!\mathbf{GGL}_S\simeq i_*\mathbf{GGL}_S$, and the result follows.
\endproof
\begin{corollary}
\label{xi_closed_sm}
Let $p:Y\to X$ be a smooth morphism and $i:Z\to Y$ be a closed immersion such that the composition $q:Z\to X$ is a closed immersion. Then $\xi_p\cdot\psi_i=\psi_q$. 
\end{corollary}

\proof
We have the following commutative diagram
\begin{equation}
\begin{split}
  \xymatrix@=10pt{
  Z \ar[r]^-{i'} \ar@{=}[rd] & Z\times_XY \ar[r]^-{q'} \ar[d]^-{p'} & Y \ar[d]^-{p}\\
  &  Z \ar[r]^-{q} & X
  }
\end{split}
\end{equation}
where $q'\circ i'=i$ and the square is Cartesian.
By Proposition~\ref{psi_comp}, Lemma~\ref{comp_ci_sm} and Lemma~\ref{xi_section} we have
\begin{equation}
\psi_q=\psi_q\cdot\xi_{p'}\cdot\psi_{i'}=\xi_p\cdot\psi_{q'}\cdot\psi_{i'}=\xi_p\cdot\psi_i
\end{equation}
and the result follows.
\endproof

\begin{corollary}
\label{xi_ci_sm}
Let $p:Y\to X$ be a smooth morphism and $i:Z\to Y$ be a closed immersion such that the composition $q:Z\to X$ is smooth. Then $\xi_p\cdot\psi_i=\xi_q$. 

\end{corollary}

\proof
We have the following commutative diagram
\begin{equation}
\begin{split}
  \xymatrix@=10pt{
  Z \ar[r]^-{i'} \ar@{=}[rd] & Z\times_XY \ar[r]^-{q'} \ar[d]^-{p'} & Y \ar[d]^-{p}\\
  &  Z \ar[r]^-{q} & X
  }
\end{split}
\end{equation}
where $q'\circ i'=i$ and the square is Cartesian.
By Lemma~\ref{xi_comp}, Lemma~\ref{xi_section} and Corollary~\ref{xi_closed_sm} we have
\begin{equation}
\xi_q=\xi_q\cdot\xi_{p'}\cdot\psi_{i'}=\xi_p\cdot\xi_{q'}\cdot\xi_{i'}=\xi_p\cdot\psi_i
\end{equation}
and the result follows.
\endproof

\begin{proposition}
\label{qpsm}
For any quasi-projective smooth morphism $p$, the map $\xi_p$ agrees with $\psi_p$.
\end{proposition}
\proof

Since every quasi-projective morphism factors into a closed immersion in an open subscheme of a projective space, by Proposition~\ref{psi_comp}, Lemma~\ref{xi_comp} and Corollary~\ref{xi_ci_sm}, it suffices to check the case where $p$ is an open immersion or the projection of a projective space. The first case follows from the definition, and the second case is already proved in Proposition~\ref{psi_proj}, which proves the result.
\endproof

\subsubsection{}
Recall that a morphism of schemes is \textbf{smoothable} if it factors as a closed immersion followed by a smooth morphism. The following proposition says that we can generalize the map $\xi_p$ to smoothable morphisms by gluing the case of closed immersions and smooth morphisms:
\begin{proposition}
\label{ext_sm_emb}
There is a unique family of isomorphisms $\xi_f:\mathbf{GGL}_Y\xrightarrow{\sim} f^!\mathbf{GGL}_X$ associated to all smoothable morphisms $f:Y\to X$ such that if $f$ factors as a closed immersion followed by a smooth morphism $Y\xrightarrow{i}Z\xrightarrow{g}X$, then $\xi_f=\xi_g\cdot\psi_i$. The map $\xi_f$ is well-defined and does not depend on the choice of factorization.
\end{proposition}
\proof
If $f:Y\to X$ factors as $Y\xrightarrow{i}Z\xrightarrow{g}X$ as stated, we set $\xi_f=\xi_g\cdot\psi_i$. If $Y\xrightarrow{i'}Z'\xrightarrow{g'}X$ is another such factorization, we have the following commutative diagram
\begin{equation}
\begin{split}
  \xymatrix@=10pt{
  Y \ar[r]^-{i_1} & Z\times_XZ' \ar[r]^-{g_1} \ar[d]^-{g_1'} & Z' \ar[d]^-{g'}\\
  &  Z \ar[r]^-{g} & X
  }
\end{split}
\end{equation}
with $i=g_1'\circ i_1$ and $i'=g_1\circ i_1$. By Lemma~\ref{xi_comp} and Corollary~\ref{xi_closed_sm} we have
\begin{equation}
\xi_g\cdot\psi_i
=
\xi_g\cdot\xi_{g_1'}\cdot\psi_{i_1}
=
\xi_{g'}\cdot\xi_{g_1}\cdot\psi_{i_1}
=
\xi_{g'}\cdot\psi_{i'}
\end{equation}
and therefore the map is well-defined and independent of the factorization. The uniqueness follows from definition.
\endproof

\subsubsection{}
Since every quasi-projective morphism factors as a closed immersion followed by a quasi-projective smooth morphism, the following result follows from Proposition~\ref{psi_comp}, Proposition~\ref{qpsm} and Proposition~\ref{ext_sm_emb}:
\begin{corollary}
\label{qproj}
For any quasi-projective morphism $p$, the map $\xi_p$ agrees with $\psi_p$.
\end{corollary}

\subsubsection{}
We are now ready to prove the general case:
\begin{proposition}
\label{sm_embed}
For any smoothable morphism $p$, the map $\xi_p$ agrees with $\psi_p$. In particular, $\xi_p=\psi_p$ for every smooth morphism $p$.
\end{proposition}
\proof

If $p:X\to S$ is a smoothable morphism, there is a non-empty open subscheme $U$ of $X$ which is quasi-projective over $S$, and the reduced complement $Z=X-U$ is smoothable over $S$. Denoting by $j:U\to X$ and $i:Z\to X$ the two immersions, since the pair of functors $(i^!,j^!)$ is conservative, the result holds for the morphism $X\to S$ if and only if it holds for both morphisms $U\to S$ and $Z\to S$. We conclude by Corollary~\ref{qproj} and noetherian induction.
\endproof

\begin{remark}
The ring spectrum $\mathbf{KGL}$ is indeed \emph{orientable} (\cite[Example 1.1.2 (4)]{Deg2}), and a fortiori so is $\mathbf{GGL}$ as a module over $\mathbf{KGL}$. But this fact has never been used throughout Chapter~\ref{chap_func}, and in particular we know that the functoriality of $\mathbf{GGL}$ from Proposition~\ref{exceptional_G} in fact does not depend on the choice of an orientation.
\end{remark}

\subsection{$G$-theory as a Borel-Moore theory}
\label{appl}
In this section we identify algebraic $G$-theory as the Borel-Moore theory associated to algebraic $K$-theory, in a way compatible with different functorialities. We fix $S\in Sch^{}_{S_0}$ a regular scheme.

\begin{recall}
\label{recall_BM}
If $\mathbb{E}$ is an oriented absolute ring spectrum (\cite[Definition 2.1.1]{Deg2}), we can associate a Borel-Moore theory as follows: for any separated morphism of finite type $p:X\to S$ and $n,m\in\mathbb{Z}$, we define
\begin{equation}
\mathbb{E}^{BM}_{n,m}(X/S):=[p_!\mathbb{S}_X(m)[n],\mathbb{E}_S]_{\mathbf{SH}(S)}=[\mathbb{S}_X(m)[n],p^!\mathbb{E}_S]_{\mathbf{SH}(X)}.
\end{equation}
It has the following natural functorialities:
\begin{itemize}
\item For every proper morphism $p:Y\to X$, we have $p_*:\mathbb{E}^{BM}_{n,m}(Y/S)\to\mathbb{E}^{BM}_{n,m}(X/S)$ given by
\begin{equation}
\begin{split}
\label{BM_comp_proper}
&\mathbb{E}^{BM}_{n,m}(Y/S)
=
[\mathbb{S}_Y(m)[n], p_Y^!\mathbb{E}_S]_{\mathbf{SH}(Y)}
\simeq
[p^*\mathbb{S}_X(m)[n], p^!p_X^!\mathbb{E}_S]_{\mathbf{SH}(X)}\\
=
&[\mathbb{S}_X(m)[n], p_*p^!p_X^!\mathbb{E}_S]_{\mathbf{SH}(X)}
\xrightarrow{ad'_{(p_!,p^!)}}
[\mathbb{S}_X(m)[n], p_X^!\mathbb{E}_S]_{\mathbf{SH}(X)}
=
\mathbb{E}^{BM}_{n,m}(X/S),
\end{split}
\end{equation}
where $p_X:X\to S$, $p_Y:Y\to S$ are structure morphisms.

\item For every smooth morphism $f:Y\to X$ of relative dimension $d$, we have $f^*:\mathbb{E}^{BM}_{n,m}(X/S)\to\mathbb{E}^{BM}_{n+2d,m+d}(Y/S)$ given by
\begin{equation}
\begin{split}
\label{BM_comp_smooth}
&\mathbb{E}^{BM}_{n,m}(X/S)
\simeq
[\mathbb{S}_X[n], p_X^!\mathbb{E}_S]_{\mathbf{SH}(X)}
\xrightarrow{ad'_{(p_\#,p^*)}}
[p_\#p^*\mathbb{S}_X[n], p_X^!\mathbb{E}_S]_{\mathbf{SH}(X)}\\
=
&[p^*\mathbb{S}_X[n], p^*p_X^!\mathbb{E}_S]_{\mathbf{SH}(Y)}
\simeq
[\mathbb{S}_Y[n], p_Y^!\mathbb{E}_S(-d)[-2d]]_{\mathbf{SH}(Y)}
=
\mathbb{E}^{BM}_{n+2d,m+d}(Y/S),
\end{split}
\end{equation}
where $p_X:X\to S$, $p_Y:Y\to S$ are structure morphisms.

\item For every regular closed immersion morphism $i:Z\to X$ of codimension $c$, we have $i^*:\mathbb{E}^{BM}_{n,m}(X/S)\to\mathbb{E}^{BM}_{n-2c,m-c}(Y/S)$ (\cite[Proposition 4.1.3]{DJK}).
\end{itemize}
\end{recall}

\subsubsection{}
In what follows, we study these maps for the ring spectrum $\mathbb{E}=\mathbf{KGL}$. The first step is to identify its Borel-Moore theory groups:
\begin{lemma}

For any separated morphism of finite type $p:X\to S$, we have an isomorphism
\begin{equation}
\label{KBM=G}
\mathbf{KGL}^{BM}_{n,m}(X/S)\simeq G_{n-2m}(X).
\end{equation}

\end{lemma}
\proof
Since $S$ is regular, by Corollary~\ref{rep_SH} and Recall~\ref{G-property} \eqref{agree_regular} there is a canonical identification $\mathbf{KGL}_S\simeq \mathbf{GGL}_S$ in $\mathbf{SH}(S)$. By Proposition~\ref{exceptional_G}, Corollary~\ref{rep_SH} and Corollary~\ref{Bott_per}, we have the following isomorphism:
\begin{align}
\begin{split}
&\mathbf{KGL}^{BM}_{n,m}(X/S)
=
[\mathbb{S}_X(m)[n],p^!\mathbf{KGL}_S]_{\mathbf{SH}(X)}
\simeq
[\mathbb{S}_X(m)[n],p^!\mathbf{GGL}_S]_{\mathbf{SH}(X)}\\
\overset{(\psi_p)^{-1}}{\simeq}
&[\mathbb{S}_X(m)[n],\mathbf{GGL}_X]_{\mathbf{SH}(X)}
\overset{\eqref{Bott_per_map}}{\simeq}
[\mathbb{S}_X[n-2m],\mathbf{GGL}_X]_{\mathbf{SH}(X)}\\
\overset{\eqref{GGL_G}}{\simeq}
&G_{n-2m}(X).
\end{split}
\end{align}
\endproof

\begin{lemma}
For every proper morphism $p:Y\to X$, the isomorphism~\eqref{KBM=G} identifies the map $f_*:\mathbf{KGL}^{BM}_{n+2m,m}(Y/S)\to\mathbf{KGL}^{BM}_{n+2m,m}(X/S)$ with the proper functoriality of $G$-theory $p_*:G_n(Y)\to G_n(X)$ (Recall~\ref{G-property} \eqref{G_proper_push}).
\end{lemma}

\proof
By the construction in Lemma~\ref{pf_GGL}, the map $\phi_p:p_*\mathbf{GGL}_Y\to \mathbf{GGL}_X$ is such that the map
\begin{equation}
\begin{split}
\label{G_GGL_prop}
G_n(Y)&\simeq[\mathbb{S}_Y[n], \mathbf{GGL}_Y]_{\mathbf{SH}(Y)}
=
[p^*\mathbb{S}_X[n], \mathbf{GGL}_Y]_{\mathbf{SH}(Y)}\\
&=[\mathbb{S}_X[n], p_*\mathbf{GGL}_Y]_{\mathbf{SH}(X)}
\xrightarrow{\phi_p}
[\mathbb{S}_X[n], \mathbf{GGL}_X]_{\mathbf{SH}(X)}
\simeq G_n(X)
\end{split}
\end{equation}
agrees with the proper functoriality of $G$-theory. Since the map $\psi_p$ is obtained from $\phi_p$ by adjunction, the map~\ref{G_GGL_prop} is induced by the proper functoriality of Borel-Moore theory~\eqref{BM_comp_proper}, and the result follows.
\endproof

\begin{lemma}
For every smooth morphism $f:Y\to X$ of relative dimension $d$, the isomorphism~\eqref{KBM=G} identifies the map $f^*:\mathbf{KGL}^{BM}_{n+2m,m}(X/S)\xrightarrow{}\mathbf{KGL}^{BM}_{n+2m+2d,m+d}(X/S)$ with the contravariant functoriality of $G$-theory $f^*:G_n(X)\xrightarrow{}G_n(Y)$ (Recall~\ref{G-property} \eqref{G_pullback}).

\end{lemma}
\proof
By construction in Definition~\ref{G_unstable} and Lemma~\ref{rep_H}, the map
\begin{equation}
G_n(X)
\simeq
[S^n\wedge X_+,G_X]_{\mathbf{H}_\bullet(X)}
\xrightarrow{}
[S^n\wedge Y_+,G_X]_{\mathbf{H}_\bullet(X)}
\simeq
G_n(Y)
\end{equation}
induced by $f_\#:Y_+\xrightarrow{}X_+$ equals the smooth functoriality of $G$-theory. By Corollary~\ref{stabilization_iso} and stabilization, we know that the map
\begin{equation}
\begin{split}
\label{G_GGL_sm}
G_n(X)
&\simeq
[\mathbb{S}_X[n], \mathbf{GGL}_X]_{\mathbf{SH}(X)}
\xrightarrow{ad'_{(p_\#,p^*)}}
[p_\#p^*\mathbb{S}_X[n], \mathbf{GGL}_X]_{\mathbf{SH}(X)}\\
&\simeq
[\Sigma^\infty Y[n], \mathbf{GGL}_X]_{\mathbf{SH}(Y)}
\simeq G_n(Y)
\end{split}
\end{equation}
agrees with the smooth functoriality of $G$-theory. It follows from Proposition~\ref{sm_embed} that the map~\eqref{G_GGL_sm} is induced by the smooth functoriality of Borel-Moore theory~\eqref{BM_comp_smooth}, and the result follows.
\endproof

\begin{lemma}
For every regular closed immersion $i:Z\to X$ of codimension $c$, the isomorphism~\eqref{KBM=G} identifies the map $i^*:\mathbf{KGL}^{BM}_{n+2m,m}(X/S)\xrightarrow{}\mathbf{KGL}^{BM}_{n+2m-2c,m-c}(X/S)$ with the contravariant functoriality of $G$-theory $i^*:G_n(X)\xrightarrow{}G_n(Z)$ (Recall~\ref{G-property} \eqref{G_pullback}).

\end{lemma}

\proof

By the construction in \cite[Definition 3.2.3]{DJK}, the Gysin morphism is given by the composition
\begin{equation}
G_n(X)
\xrightarrow{\gamma_t}
G_{n+1}(\mathbb{G}_{m,X})
\xrightarrow{\partial}
G_n(N_ZX)
\simeq
G_n(Z).
\end{equation}
where $\gamma_t:G_n(X)\to G_{n+1}(\mathbb{G}_{m,X})$ is given by the multiplication by $t\in K_1(\mathbb{Z}[t,t^{-1}])$, and $\partial$ is the boundary map in the long exact sequence of homotopy groups associated to the homotopy fiber sequence $G(N_ZX)\to G(D_ZX)\to G(\mathbb{G}_{m,X})$, $D_ZX=Bl_{\mathbb{A}^1_Z}\mathbb{A}^1_X-Bl_ZX$ being the deformation space to the normal cone.
\footnote{Equivalently, under the canonical isomorphism $G_{n+1}(\mathbb{G}_{m,X})\simeq G_n(X)\oplus G_{n+1}(X)$, $\gamma_t$ is the inclusion of the first summand, see \cite[V.6.2]{Wei}.}
We claim that the composition $G_n(X)
\xrightarrow{\gamma_t}
G_{n+1}(\mathbb{G}_{m,X})
\xrightarrow{\partial}
G_n(N_ZX)$ is given by the contravariant functoriality $G_n(X)\to G_n(N_ZX)$. Indeed, the very argument of \cite[Proposition 2.30]{Jin} applies: using the (double) deformation to the normal cone, we are reduced to the case where $i:Z\to\mathbb{A}^1_Z$ is the zero section, and in this case $\partial$ is a retraction of $\gamma_t$.
\footnote{Alternatively, we can apply the same argument by looking at Borel-Moore motives in the homotopy category of $\mathbf{KGL}$-modules.}
The result then follows from the naturality of the contravariant functoriality of $G$-theory.
\endproof

\begin{corollary}
\label{G=BM_func}
For any separated scheme of finite type $X$ over $S$, there is a natural isomorphism 
\begin{equation}
\mathbf{KGL}^{BM}_{n,m}(X/S)\simeq G_{n-2m}(X)
\end{equation}
such that
\begin{itemize}
\item The proper functoriality on $\mathbf{KGL}^{BM}$ for proper morphisms is given by the proper functoriality of $G$-theory;
\item The Gysin morphisms on $\mathbf{KGL}^{BM}$ for lci morphisms (\cite[Definition 3.3.2]{Deg2}) is given by the contravariant functoriality of $G$-theory.
\end{itemize}

\end{corollary}

\begin{remark}
\label{remark_RR}
Corollary~\ref{G=BM_func} has been used in \cite[Example 3.3.11 (2)]{Deg2} to obtain a Riemann-Roch theorem for singular schemes using $G$-theory: by formal properties of Borel-Moore theories, one constructs a map of spectra 
\begin{equation}
\mathrm{ch}:\mathbf{GGL}\to\oplus_{i\in\mathbb{Z}}\mathbf{H}\mathbb{Q}^{BM}(i)[2i]
\end{equation}
where the right hand side is the spectrum representing Borel-Moore motivic homology with rational coefficients, which is a Chern character map that generalizes the Chern character at the level of motivic spectra in \cite[Definition 6.2.3.9]{Rio2} to singular schemes. Using the construction of Gysin morphisms (\cite[Definition 3.3.2]{Deg2}, which has been extended to a more general setting in \cite{DJK}), it is proved that the map $\mathrm{ch}$ satisfies compatibilities with proper functoriality and Gysin morphisms, which generalizes the statement in \cite[Theorem 18.3]{Ful}. This result is proved mostly with tools from $\mathbb{A}^1$-homotopy theory and the six functors formalism, without using MacPherson's ``graph construction'' method in \cite[\S18]{Ful}.

\end{remark}

\end{document}